\documentclass[11pt]{article}

\usepackage{amsmath}
\usepackage{amsthm}
\usepackage{latexsym}
\usepackage{amssymb}
\usepackage{graphicx}
\usepackage{amsmath}
\usepackage{amssymb}
\usepackage{graphicx}
\usepackage{amsmath}
\usepackage{amsfonts}
\usepackage{euscript}
\usepackage{epsfig}
 \usepackage{color}
\usepackage{amsmath}
\usepackage{latexsym}
\usepackage{amssymb}
\usepackage{cite}
\usepackage{indentfirst}
\usepackage{enumerate}
\usepackage{tikz}

\newtheorem{thm}{Theorem}
\newtheorem{lem}{Lemma}
\newtheorem{cor}{Corollary}

\title{\bf Invariants
 of Third Type Almost Geodesic Mappings\\
 of Generalized Riemannian Space}
\author{Nenad O. Vesi\'c}
\date{}

\numberwithin{equation}{section}

\makeatletter
\def\maketag@@@#1{\hbox{\m@th\normalfont\normalsize#1}}
\makeatother
\usepackage{lipsum}

\newcommand\blfootnote[1]{%
  \begingroup
  \renewcommand\thefootnote{}\footnote{#1}%
  \addtocounter{footnote}{-1}%
  \endgroup
}
\usepackage{empheq}

\begin{document}
  \maketitle

  \begin{abstract}
    We studied rules of transformations of
    Christoffel symbols under third type almost geodesic mappings
     in this paper. From this research, we obtained some new
     invariants of these mappings. These invariants are analogies of
     Thomas projective parameter and Weyl projective tensor.
     \vspace{.2cm}


    \textbf{$2010$ Math. Subj. Classification:} 53C15, 47A15, 58C30, 55C99, 53A55, 35R01
    \blfootnote{This paper is financially supported by Serbian
    Ministry of Education, Science and Technological Development,
    Grant No. 174012}
  \end{abstract}

  \section{Introduction}

  Following the Eisenhart's work \cite{eis1, eis01, eis02}, it is
  started the research about Riemannian spaces endowed with
  non-symmetric metrics \cite{mileva1, mileva2}.
  An $N$-dimensional manifold endowed with metric tensor $g_{ij}$
  non-symmetric by indices $i$ and $j$ is \textbf{the generalized
  Riemannian space $\mathbb{GR}_N$}. Affine connection
  coefficients of the space $\mathbb{GR}_N$ are Christoffel
  symbols of the second kind $\Gamma^i_{jk}$ with respect to the
  connection of the metric $g_{ij}$. Because
  $\Gamma^i_{jk}\neq\Gamma^i_{kj}$, the symmetric and anti-symmetric
  parts of the Christoffel symbol $\Gamma^i_{jk}$ are defined as

  \begin{eqnarray}
  \Gamma^i_{\underline{jk}}=
  \frac12(\Gamma^i_{jk}+\Gamma^i_{kj})&\mbox{and}&
  \Gamma^i_{\underset\vee{jk}}=
  \frac12(\Gamma^i_{jk}-\Gamma^i_{kj}).
  \end{eqnarray}

  \noindent The anti-symmetric part $\Gamma^i_{\underset\vee{jk}}$ of
  the coefficient $\Gamma^i_{jk}$ is called \textbf{the torsion tensor
  of the space $\mathbb{GR}_N$}. The Riemannian space $\mathbb R_N$,
  endowed with affine connection coefficients
  $\Gamma^i_{\underline{jk}}$ is \textbf{the associated space of
  the space $\mathbb{GR}_N$}
  \cite{marijamilanmica,micaljubicamincicmilan,micamilanljubica2,micamincic,
  micamincicljubica1,micamincicljubica2}.

  With regard to the affine connection of Riemannian space, it is
  defined one kind of covariant derivation

  \begin{equation}
  a^i_{j;k}=a^i_{j,k}+\Gamma^i_{\underline{\alpha k}}a^\alpha_j-
  \Gamma^\alpha_{\underline{jk}}a^i_\alpha,
  \label{eq:covderivativeAn}
\end{equation}
for a tensor $a^i_j$ of the type $(1,1)$ and partial derivation
denoted by comma. Curvature tensor $R^i_{jmn}$ of the associated
space $\mathbb R_N$ is

\begin{equation}
  R^i_{jmn}=\Gamma^i_{\underline{jm},n}-
  \Gamma^i_{\underline{jn},m}+
  \Gamma^\alpha_{\underline{jm}}\Gamma^i_{\underline{\alpha n}}-
  \Gamma^\alpha_{\underline{jn}}\Gamma^i_{\underline{\alpha m}}.
  \label{eq:r}
\end{equation}

Four kinds of covariant differentiation with regard to affine
connection of the space $\mathbb{GR}_N$ are \cite{mincabc}:

\begin{eqnarray}
  a^i_{j\underset1|k}=a^i_{j,k}+\Gamma^i_{\alpha k}a^\alpha_j-
  \Gamma^\alpha_{jk}a^i_\alpha&&
  a^i_{j\underset2|k}=a^i_{j,k}+\Gamma^i_{k\alpha}a^\alpha_j-
  \Gamma^\alpha_{kj}a^i_\alpha,\label{eq:covdev1}\\
  a^i_{j\underset3|k}=a^i_{j,k}+\Gamma^i_{\alpha k}a^\alpha_j-
  \Gamma^\alpha_{kj}a^i_\alpha&&
  a^i_{j\underset4|k}=a^i_{j,k}+\Gamma^i_{\alpha k}a^\alpha_j-
  \Gamma^\alpha_{kj}a^i_\alpha.\label{eq:covdev4}
\end{eqnarray}

\noindent There are twelve curvature tensors of the space
$\mathbb{GR}_N$ \cite{mincabc}. They are elements of the family

\begin{equation}
    K^i_{jmn}=R^i_{jmn}+u\Gamma^i_{\underset\vee{jm};n}+
  u'\Gamma^i_{\underset\vee{jn};m}+
  v\Gamma^\alpha_{\underset\vee{jm}}\Gamma^i_{\underset\vee{\alpha
  n}}+v'\Gamma^\alpha_{\underset\vee{jn}}\Gamma^i_{\underset\vee{\alpha
  m}}+w\Gamma^\alpha_{\underset\vee{mn}}\Gamma^i_{\underset\vee{\alpha
  j}},\label{eq:R}
\end{equation}

\noindent for real constants $u,u',v,v',w$. Five of these curvature
tensors are linearly independent \cite{wc10}.

  \subsection{Mappings of generalized Riemannian spaces}

  Riemannian and generalized Riemannian spaces are special affine
  connected spaces. We are going to pay attention on mappings
  between generalized Riemannian spaces in this paper.

  Different mappings of Riemannian and generalized Riemannian spaces as well
  as their invariants have been investigated (see \cite{mikesberezovski1, marijamilanmica,
   eis1,  hall1, hall2, zlativanov, mikes1, mikesnovi, mikes2, mincabc,
  wc10, mincmica, mincicmica, sinjukov, micapi1, micapi2, micapi3,
  micamincic, micamincicljubica2, micamincicljubica1, micaljubicamincicmilan,
  micamilanljubica2, jamicapi3}).
 N. S. Sinyukov started the research about almost geodesic mappings \cite{sinjukov}.
  His work has been continued by by J. Mike\v s and his research group (see \cite{mikes1, mikes2,
  mikesnovi}).

  Sinyukov \cite{sinjukov} generalized the term of geodesics. He
  involved the terms of almost geodesic lines and almost geodesic
  mappings of symmetric affine connection spaces. He  founded that
   are three types $\pi_1,\pi_2,\pi_3$ of almost geodesic of a
  symmetric affine connection space.

  The terms of almost geodesic lines and almost geodesic mappings
  are generalized in \cite{micapi1, micapi2, micapi3}. It is founded
  that there are three types and two kinds of almost geodesic
  mappings of a non-symmetric affine connection space. Because
  generalized Riemannian spaces are special non-symmetric affine
  connection spaces, the basic equations of an almost geodesic
  mapping $f:\mathbb{GR}_N\rightarrow\mathbb{G\overline R}_N$
   of the third type and $s$-th kind, $s=1,2$, are:

    \begin{align}
    &\underset s\pi{}_3:\left\{\begin{array}{l}
    \overline
    \Gamma^i_{jk}=\Gamma^i_{jk}+\psi_j\delta^i_k+\psi_k\delta^i_j+2\sigma_{jk}\varphi^i+
    \xi^i_{jk},\\
    \varphi^i_{\underset s|j}+(-1)^{s-1}\xi^i_{\alpha j}\varphi^\alpha=
    \nu_j\varphi^i+\mu\delta^i_j,
  \end{array}\right.\label{eq:pi3GRNnsimbasic}
  \end{align}

  \noindent for tensor $\xi^i_{jk}$ anti-symmetric by indices $j$
  and $k$, covariant vectors $\psi_j,\nu_j$, contravariant vector
  $\varphi^i$ and scalar function $\mu$. Third type almost
  geodesic mapping $f$ has the property of reciprocity if its
  inverse mapping is third type almost geodesic mapping.

  In this paper, we will pay attention to invariants of equitorsion
  third type almost geodesic
  mappings which satisfy the property of reciprocity.

  \subsection{Motivation}

  A. Einstein was the first scientist who applied the non-symmetric
  affine connection in the theory of gravitation \cite{e1,e2,e3}. In
  his theory, Weyl projective tensor is related to gravity.

  In this paper, we wish to generalize Weyl projective tensor as an
  invariant of equitorsion third type almost geodesic mappings.
  It will be obtained transformation rules of covariantly
  differentiated torsion tensor $\Gamma^i_{\underset\vee{jk}}$ under
  equitorsion third type almost geodesic mappings
  which have the property of reciprocity.
  In the next, we will find families of
  invariants of these mappings with regard
  to the changes of curvature tensors $K^i_{jmn}$ given by the
  equation (\ref{eq:R}).

  \section{Main results}

  Let $f:\mathbb{GR}_N\rightarrow\mathbb{G\overline R}_N$ be an
  equitorsion almost geodesic mapping of the third type and $s$-th kind $s=1,2$, which
  satisfies the property of reciprocity. Basic equations of this
  mapping are

  \begin{align}
    &\underset s\pi{}_3:\left\{\begin{array}{l}
    \overline
    \Gamma^i_{jk}=\Gamma^i_{jk}+\psi_j\delta^i_k+\psi_k\delta^i_j+2\sigma_{jk}\varphi^i,\\
    \varphi^i_{\underset s|j}=\nu_j\varphi^i+\mu\delta^i_j,
  \end{array}\right.\label{eq:pi3nsimbasic}
  \end{align}

  \noindent for scalar function $\mu$, covariant vectors $\psi_j,\nu_j$,
  contravariant vector $\varphi^i$ and tensor $\sigma_{jk}$ of the type
  $(0,2)$ symmetric by indices $j$ and $k$. It is obtained in
  \cite{jaPI3} that the geometrical objects

  \begin{align}
    &\aligned
    \underset{(1)}{\overset\star{\mathcal
    W}}{}^i_{jmn}&=R^i_{jmn}+\delta^i_j\underset{(1)}{\overset\star\eta}{}_{[mn]}
    -\frac1{N+1}\delta^i_m\big(
    \Gamma^\alpha_{\underline{j\alpha};n}-(N+1)(\underset{(1)}{\overset\star\eta}{}_{jn}+
    \mu\sigma_{jn})\big)\\&+
    \frac1{N+1}\delta^i_n\big(\Gamma^\alpha_{\underline{j\alpha};m}-
    (N+1)(\underset{(1)}{\overset\star\eta}{}_{jm}+\mu\sigma_{jm})\big)\\&
    -
    \big(\sigma_{jm;n}-\sigma_{jn;m}-(\sigma_{jm}\sigma_{\alpha n}-
    \sigma_{jn}\sigma_{\alpha m})\varphi^\alpha\big)\varphi^i+\sigma_{jm}\Gamma^i_{\underset\vee{\alpha
    n}}\varphi^\alpha-
    \sigma_{jn}\Gamma^i_{\underset\vee{\alpha m}}\varphi^\alpha,
    \endaligned\label{eq:inv1}\\
    &\aligned
    \underset{(2)}{\overset\star{\mathcal
    W}}{}^i_{jmn}&=R^i_{jmn}+\delta^i_j\underset{(2)}{\overset\star\eta}{}_{[mn]}-\frac1{N+1}\delta^i_m\big(
    \Gamma^\alpha_{\underline{j\alpha};n}-(N+1)(\underset{(2)}{\overset\star\eta}{}_{jn}+
    \mu\sigma_{jn})\big)\\&+
    \frac1{N+1}\delta^i_n\big(\Gamma^\alpha_{\underline{j\alpha};m}-
    (N+1)(\underset{(2)}{\overset\star\eta}{}_{jm}+\mu\sigma_{jm})\big)\\&
    -
    \big(\sigma_{jm;n}-\sigma_{jn;m}-(\sigma_{jm}\sigma_{\alpha n}-
    \sigma_{jn}\sigma_{\alpha m})\varphi^\alpha\big)\varphi^i-\sigma_{jm}\Gamma^i_{\underset\vee{\alpha
    n}}\varphi^\alpha+
    \sigma_{jn}\Gamma^i_{\underset\vee{\alpha m}}\varphi^\alpha,
    \endaligned\label{eq:inv2}
  \end{align}

  \noindent for $\underset{(1)}{\overset\star\eta}{}_{ij}$ and $\underset{(2)}{\overset\star\eta}{}_{ij}$
    given by the following equations:

    \begin{align}
  &\aligned
  \underset{(1)}{\overset\star\eta}{}_{jk}&=\frac1{(N+1)^2}\Big((N+1)\varphi^\alpha\sigma_{jk}
  (\Gamma^\beta_{\underline{\alpha\beta}}+
  \sigma_{\alpha\beta}\varphi^\beta)-(\Gamma^\alpha_{\underline{j\alpha}}+
  \sigma_{j\alpha}\varphi^\alpha)(\Gamma^\beta_{\underline{k\beta}}+\sigma_{k\beta}\varphi^\beta)\Big)\\&-
  \frac1{N+1}\big(\sigma_{j\alpha;k}\varphi^\alpha+\sigma_{jk}\mu+\sigma_{j\alpha}(
  \nu_k\varphi^\alpha-\Gamma^\alpha_{\underset\vee{\beta k}}\varphi^\beta)\big),
  \endaligned\label{eq:eta1}\\
  &\aligned
  \underset{(2)}{\overset\star\eta}{}_{jk}&=\frac1{(N+1)^2}\Big((N+1)\varphi^\alpha\sigma_{jk}
  (\Gamma^\beta_{\underline{\alpha\beta}}+
  \sigma_{\alpha\beta}\varphi^\beta)-(\Gamma^\alpha_{\underline{j\alpha}}+
  \sigma_{j\alpha}\varphi^\alpha)(\Gamma^\beta_{\underline{k\beta}}+\sigma_{k\beta}\varphi^\beta)\Big)\\&-
  \frac1{N+1}\big(\sigma_{j\alpha;k}\varphi^\alpha+\sigma_{jk}\mu+\sigma_{j\alpha}(
  \nu_k\varphi^\alpha+\Gamma^\alpha_{\underset\vee{\beta
  k}}\varphi^\beta)\big),
  \endaligned\label{eq:eta2}
\end{align}

\noindent     scalar function $\mu$ and antisymmetrization
    without division denoted by square brackets, are invariants of the mapping
    $f$.

  \subsection{Transformations of covariant derivative of torsion tensor}

  Let $f:\mathbb{GR}_N\rightarrow\mathbb{G\overline R}_N$ be an
  equitorsion third type almost geodesic mapping of an $s$-th kind,
  $s=1,2$,
  which has the property of reciprocity.
  Based on the equation (\ref{eq:covderivativeAn}) and the
  invariance
  $\overline\Gamma^i_{\underset\vee{jk}}=\Gamma^i_{\underset\vee{jk}}$,
  we obtain that is

  \begin{equation}
    \aligned
    \overline\Gamma^i_{\underset\vee{jm}\overline;n}-
    \Gamma^i_{\underset\vee{jm};n}&=\overline\Gamma^i_{\underline{\alpha
    n}}\overline\Gamma^\alpha_{\underset\vee{jm}}-
    \overline\Gamma^\alpha_{\underline{jn}}\overline\Gamma^i_{\underset\vee{\alpha
    m}}-
    \overline\Gamma^\alpha_{\underline{mn}}\overline\Gamma^i_{\underset\vee{j\alpha}}-
    \Gamma^i_{\underline{\alpha
    n}}\Gamma^\alpha_{\underset\vee{jm}}+
    \Gamma^\alpha_{\underline{jn}}\Gamma^i_{\underset\vee{\alpha
    m}}+
    \Gamma^\alpha_{\underline{mn}}\Gamma^i_{\underset\vee{j\alpha}}\\&=
    \Gamma^\alpha_{\underset\vee{jm}}(\overline\Gamma^i_{\underline{\alpha n}}-
    \Gamma^i_{\underline{\alpha n}})-\Gamma^i_{\underset\vee{\alpha
    m}}(\overline\Gamma^\alpha_{\underline{jn}}-\Gamma^\alpha_{\underline{jn}})-
    \Gamma^i_{\underset\vee{j\alpha}}(\overline\Gamma^\alpha_{\underline{mn}}-
    \Gamma^\alpha_{\underline{mn}})
    \endaligned\label{eq:T-Tpi3}
  \end{equation}

  Because the mapping $f$ has the property of reciprocity, it is
  got in \cite{jaPI3} that is

  \begin{eqnarray}
    \overline\psi_j=-\psi_j&\mbox{and}&
    \overline\sigma_{jk}\overline\varphi^i=-
    \sigma_{jk}\varphi^i,
  \end{eqnarray}

  \noindent i.e.

  \begin{equation}
    \aligned
    \overline\Gamma^i_{\underline{jk}}-
    \Gamma^i_{\underline{jk}}&=\frac1{N+1}\big((\overline\Gamma^\alpha_{\underline{j\alpha}}+\overline\sigma_{j\alpha}\overline\varphi^\alpha)\delta^i_k+
    (\overline\Gamma^\alpha_{\underline{k\alpha}}+\overline\sigma_{k\alpha}\overline\varphi^\alpha)\delta^i_j\big)-
    \overline\sigma_{jk}\overline\varphi^i
    \\&-
    \frac1{N+1}\big((\Gamma^\alpha_{\underline{j\alpha}}+\sigma_{j\alpha}\varphi^\alpha)\delta^i_k+
    (\Gamma^\alpha_{\underline{k\alpha}}+\sigma_{k\alpha}\varphi^\alpha)\delta^i_j\big)+
    \sigma_{jk}\varphi^i.
    \endaligned\label{eq:pi3L-Lfactorized}
  \end{equation}

  With regard to the expressions (\ref{eq:T-Tpi3},
  \ref{eq:pi3L-Lfactorized}), we obtain that it holds

  \begin{equation}
    \overline\Gamma^i_{\underset\vee{jm}\overline;n}-
    \Gamma^i_{\underset\vee{jm};n}=\underset{(p)}{\overline\sigma}{}^i_{jmn}-
    \underset{(p)}\sigma{}^i_{jmn},
    \label{eq:Gamma-Gamma=sigma-sigmapi3}
  \end{equation}

  \noindent $p=1,\ldots,8$, for

  \begin{align}
    &\aligned
    \underset{(1)}\sigma{}^i_{jmn}&=\Gamma^\alpha_{\underset\vee{jm}}
    \Gamma^i_{\underline{\alpha n}}-\Gamma^i_{\underset\vee{\alpha
    m}}\Gamma^\alpha_{\underline{jn}}-\Gamma^i_{\underset\vee{j\alpha}}
    \Gamma^\alpha_{\underline{mn}},
    \endaligned\label{eq:pi3sigma1}\\\displaybreak[0]
    &\aligned
    \underset{(2)}\sigma{}^i_{jmn}&=\Gamma^\alpha_{\underset\vee{jm}}
    \Gamma^i_{\underline{\alpha n}}+\Gamma^i_{\underset\vee{\alpha
    m}}\varphi^\alpha\sigma_{jn}+\Gamma^i_{\underset\vee{j\alpha}}\varphi^\alpha
    \sigma_{mn}-\frac1{N+1}\big(
    2\Gamma^i_{\underset\vee{jm}}\Gamma^\alpha_{\underline{n\alpha}}+
    \Gamma^i_{\underset\vee{jn}}\Gamma^\alpha_{\underline{m\alpha}}-
    \Gamma^i_{\underset\vee{mn}}\Gamma^\alpha_{\underline{j\alpha}}
    \big)\\&-
    \frac1{N+1}\big(
    2\Gamma^i_{\underset\vee{jm}}\varphi^\alpha\sigma_{\alpha n}+
    \Gamma^i_{\underset\vee{jn}}\varphi^\alpha\sigma_{\alpha m}-
    \Gamma^i_{\underset\vee{mn}}\varphi^\alpha\sigma_{\alpha j}
    \big),
    \endaligned\label{eq:pi3sigma2}\\\displaybreak[0]
    &\aligned
    \underset{(3)}\sigma{}^i_{jmn}&=\Gamma^\alpha_{\underset\vee{jm}}
    \Gamma^i_{\underline{\alpha n}}-\Gamma^i_{\underset\vee{\alpha
    m}}\Gamma^\alpha_{\underline{jn}}+
    \Gamma^i_{\underset\vee{j\alpha}}\varphi^\alpha\sigma_{mn}\\&-
    \frac1{N+1}\big(\Gamma^i_{\underset\vee{jm}}
    \Gamma^\alpha_{\underline{n\alpha}}+
    \Gamma^i_{\underset\vee{jn}}\Gamma^\alpha_{\underline{m\alpha}}
    +\Gamma^i_{\underset\vee{jm}}\varphi^\alpha\sigma_{\alpha n}+
    \Gamma^i_{\underset\vee{jn}}\varphi^\alpha\sigma_{\alpha
    m}\big),
    \endaligned\label{eq:pi3sigma3}\\\displaybreak[0]
    &\aligned
    \underset{(4)}\sigma{}^i_{jmn}&=\Gamma^\alpha_{\underset\vee{jm}}
    \Gamma^i_{\underline{\alpha
    n}}-\Gamma^i_{\underset\vee{j\alpha}}\Gamma^\alpha_{\underline{mn}}+
    \Gamma^i_{\underset\vee{\alpha m}}\varphi^\alpha\sigma_{jn}\\&-
    \frac1{N+1}\big(\Gamma^i_{\underset\vee{jm}}\Gamma^\alpha_{\underline{n\alpha}}-
    \Gamma^i_{\underset\vee{mn}}\Gamma^\alpha_{\underline{j\alpha}}+
    \Gamma^i_{\underset\vee{jm}}\varphi^\alpha\sigma_{\alpha
    n}-\Gamma^i_{\underset\vee{mn}}\varphi^\alpha\sigma_{\alpha j}
    \big),
    \endaligned\label{eq:pi3sigma4}\\\displaybreak[0]
    &\aligned
    \underset{(5)}\sigma{}^i_{jmn}&=-\Gamma^\alpha_{\underset\vee{jm}}
    \varphi^i\sigma_{\alpha n}-\Gamma^i_{\underset\vee{\alpha
    m}}\Gamma^\alpha_{\underline{jn}}-\Gamma^i_{\underset\vee{j\alpha}}
    \Gamma^\alpha_{\underline{mn}}\\&+
    \frac1{N+1}\big(\delta^i_n\Gamma^\alpha_{\underset\vee{jm}}\Gamma^\beta_{\underline{\alpha\beta}}
    +\Gamma^i_{\underset\vee{jm}}\Gamma^\alpha_{\underline{n\alpha}}+
    \delta^i_n\Gamma^\alpha_{\underset\vee{jm}}\varphi^\beta\sigma_{\alpha
    \beta}+\Gamma^i_{\underset\vee{jm}}\varphi^\alpha\sigma_{\alpha n}
    \big),
    \endaligned\label{eq:pi3sigma5}\\\displaybreak[0]
    &\aligned
    \underset{(6)}\sigma{}^i_{jmn}&=-\Gamma^\alpha_{\underset\vee{jm}}
    \varphi^i\sigma_{\alpha n}
    +\Gamma^i_{\underset\vee{\alpha
    m}}\varphi^\alpha\sigma_{jn}+\Gamma^i_{\underset\vee{j\alpha}}\varphi^\alpha
    \sigma_{mn}\\&+\frac1{N+1}\big(\delta^i_n
    \Gamma^\alpha_{\underset\vee{jm}}\Gamma^\beta_{\underline{\alpha\beta}}-
    \Gamma^i_{\underset\vee{jm}}\Gamma^\alpha_{\underline{n\alpha}}-
    \Gamma^i_{\underset\vee{jn}}\Gamma^\alpha_{\underline{m\alpha}}+
    \Gamma^i_{\underset\vee{mn}}\Gamma^\alpha_{\underline{j\alpha}}
    \big)\\&+
    \frac1{N+1}\big(\delta^i_n\Gamma^\alpha_{\underset\vee{jm}}\varphi^\beta\sigma_{\alpha\beta}-
    \Gamma^i_{\underset\vee{jm}}\varphi^\alpha\sigma_{\alpha n}-
    \Gamma^i_{\underset\vee{jn}}\varphi^\alpha\sigma_{\alpha m}+
    \Gamma^i_{\underset\vee{mn}}\varphi^\alpha\sigma_{\alpha j}
    \big),
    \endaligned\label{eq:pi3sigma6}\\\displaybreak[0]
    &\aligned
    \underset{(7)}\sigma{}^i_{jmn}&=-\Gamma^\alpha_{\underset\vee{jm}}
    \varphi^i\sigma_{\alpha n}-\Gamma^i_{\underset\vee{\alpha
    m}}\Gamma^\alpha_{\underline{jn}}+
    \Gamma^i_{\underset\vee{j\alpha}}\varphi^\alpha\sigma_{mn}\\&+
    \frac1{N+1}\big(\delta^i_n
    \Gamma^\alpha_{\underset\vee{jm}}\Gamma^\beta_{\underline{\alpha\beta}}-
    \Gamma^i_{\underset\vee{jn}}\Gamma^\alpha_{\underline{m\alpha}}
    \big)+
    \frac1{N+1}\big(
    \delta^i_n\Gamma^\alpha_{\underset\vee{jm}}\varphi^\beta\sigma_{\alpha\beta}-
    \Gamma^i_{\underset\vee{jn}}\varphi^\alpha\sigma_{\alpha
    m}\big),
    \endaligned\label{eq:pi3sigma7}\\\displaybreak[0]
    &\aligned
    \underset{(8)}\sigma{}^i_{jmn}&=-\Gamma^\alpha_{\underset\vee{jm}}
    \varphi^i\sigma_{\alpha n}-\Gamma^i_{\underset\vee{j\alpha}}\Gamma^\alpha_{\underline{mn}}+
    \Gamma^i_{\underset\vee{\alpha m}}\varphi^\alpha\sigma_{jn}\\&+
    \frac1{N+1}\big(\delta^i_n
    \Gamma^\alpha_{\underset\vee{jm}}\Gamma^\beta_{\underline{\alpha\beta}}+
    \Gamma^i_{\underset\vee{mn}}\Gamma^\alpha_{\underline{j\alpha}}
    \big)+\frac1{N+1}\big(\delta^i_n\Gamma^\alpha_{\underset\vee{jm}}
    \varphi^\beta\sigma_{\alpha\beta}+
    \Gamma^i_{\underset\vee{mn}}\varphi^\alpha\sigma_{\alpha j}\big)
    ,
    \endaligned\label{eq:pi3sigma8}
  \end{align}

  \noindent and the corresponding
  $\underset{(p)}{\overline\sigma}{}^i_{jmn}$.

  Let be

  \begin{equation}
  \aligned
    &\underset1U=\Gamma^\alpha_{\underset\vee{jm}}
    \Gamma^i_{\underline{\alpha n}},\quad
    \underset2U=\Gamma^\alpha_{\underset\vee{jn}}
    \Gamma^i_{\underline{\alpha m}},
    \quad\underset3U=
    \Gamma^i_{\underset\vee{\alpha
    m}}\Gamma^\alpha_{\underline{jn}},\quad
    \underset4U=\Gamma^i_{\underset\vee{\alpha n}}\Gamma^\alpha_{\underline{jm}},\quad
    \underset5U=\Gamma^i_{\underset\vee{j\alpha}}\Gamma^\alpha_{\underline{mn}},\quad
    \\&\underset6U=\Gamma^i_{\underset\vee{jm}}\Gamma^\alpha_{\underline{n\alpha}},\quad
    \underset7U=\Gamma^i_{\underset\vee{jn}}\Gamma^\alpha_{\underline{m\alpha}},\quad
    \underset8U=\Gamma^i_{\underset\vee{mn}}\Gamma^\alpha_{\underline{j\alpha}},\quad
    \underset9U=\Gamma^i_{\underset\vee{\alpha m}}\varphi^\alpha\sigma_{jn},\quad
    \underset{10}U=\Gamma^i_{\underset\vee{\alpha
    n}}\varphi^\alpha\sigma_{jm},\\&
    \underset{11}U=\Gamma^i_{\underset\vee{j\alpha}}\varphi^\alpha\sigma_{mn},\quad
    \underset{12}U=\Gamma^i_{\underset\vee{jm}}\varphi^\alpha\sigma_{\alpha
    n},\quad
    \underset{13}U=\Gamma^i_{\underset\vee{jn}}\varphi^\alpha\sigma_{\alpha
    m},\quad
    \underset{14}U=\Gamma^i_{\underset\vee{mn}}\varphi^\alpha\sigma_{\alpha
    j},\\&
    \underset{15}U=\delta^i_n\Gamma^\alpha_{\underset\vee{jm}}
    \Gamma^\beta_{\underline{\alpha\beta}},\quad
    \underset{16}U=\delta^i_m\Gamma^\alpha_{\underset\vee{jn}}
    \Gamma^\beta_{\underline{\alpha\beta}},
    \underset{17}U=\delta^i_n\Gamma^\alpha_{\underset\vee{jm}}\varphi^\beta\sigma_{\alpha\beta},\quad
    \underset{18}U=\delta^i_m\Gamma^\alpha_{\underset\vee{jn}}\varphi^\beta\sigma_{\alpha\beta},\\&
    \underset{19}U=\Gamma^\alpha_{\underset\vee{jm}}\varphi^i\sigma_{\alpha n},\quad
    \underset{20}U=\Gamma^\alpha_{\underset\vee{jn}}\varphi^i\sigma_{\alpha m}.
  \endaligned\label{eq:UthetaPI3}
  \end{equation}

  \noindent It holds the following lemma:

  \begin{lem}
    Let $f:\mathbb{GR}_N\rightarrow\mathbb{G\overline R}_N$ be an
    equitorsion third type almost geodesic mapping which has the
    property of reciprocity.

    \begin{enumerate}[\emph{aa)}]
      \item[\emph{a)}] Covariant derivatives
      $\Gamma^i_{\underset\vee{jm};n}$ and
      $\overline\Gamma^i_{\underset\vee{jm}\overline;n}$ of the torsion tensor of the
      spaces $\mathbb{GR}_N$ and $\mathbb{G\overline R}_N$ satisfy the
      equations

      \begin{equation}
        \overline\Gamma^i_{\underset\vee{jm}\overline;n}=
        \Gamma^i_{\underset\vee{jm};n}+\underset{(p)}{\overline\sigma}{}^i_{jmn}-
        \underset{(p)}\sigma{}^i_{jmn}=
        \Gamma^i_{\underset\vee{jm};n}+\sum_{\rho=1}^8{
        \sum_{\theta=1}^{20}{u^\rho_{\theta}\big(\underset\theta{\overline U}-
        \underset\theta U\big)}},
        \label{eq:T;toT;pi3}
      \end{equation}

      \noindent $p=1,\ldots,8$, for the corresponding
      real constants $u^\rho_\theta$, geometrical
      objects $\underset{(p)}\sigma{}^i_{jmn},\underset\theta U$
      given by the equations
      \emph{(\ref{eq:pi3sigma1}--\ref{eq:UthetaPI3})}.
      \item[\emph{b)}] The rank of matrix $\big[u^\rho_{\theta}\big]_{8\times20},
      \rho=1,\ldots,8$, is $4$, i.e. there are four linearly independent
      transformations of the transformations from
      \emph{(\ref{eq:T;toT;pi3})}.\qed
    \end{enumerate}
  \end{lem}

  \begin{cor}
    Geometrical objects

    \begin{equation}
      \overset{\rho}{\widetilde{\mathcal
      T}}{}^i_{jm;n}=\Gamma^i_{\underset\vee{jm};n}-{
        \sum_{\theta=1}^{20}{u^{\rho}_{\theta}
        \underset\theta U}}{}^i_{jmn},
      \label{eq:T;inv}
    \end{equation}

    \noindent $\rho\in\{1,\ldots,8\}$, for the corresponding real constants $u^\rho_\theta$, are invariants of an equitorsion almost geodesic
    mapping $f:\mathbb{GR}_N\rightarrow\mathbb{G\overline R}_N$
    which has the property of reciprocity. Four of these invariants are linearly independent. \qed
  \end{cor}




  \subsection{Transformations of curvature tensors under almost geodesic mappings}

  Let $f:\mathbb{GR}_N\rightarrow\mathbb{G\overline R}_N$ be an
  equitorsion almost geodesic mapping of the third type and $s$-th
  kind, $s=1,2$, which has the property of reciprocity.
  From the invariance of the geometrical objects
  $\underset{(1)}{\overset\ast{\mathcal W}}{}^i_{jmn}$ and
  $\underset{(2)}{\overset\ast{\mathcal W}}{}^i_{jmn}$ given by the
  equations (\ref{eq:inv1}, \ref{eq:inv2}), we directly
  obtain that is

  \begin{align}
    &\aligned
    \overline R^i_{jmn}&=R^i_{jmn}-
    \delta^i_j\underset{(1)}{\overset\star{\overline\eta}}{}_{[mn]}
    +\frac1{N+1}\delta^i_m\big(
    \overline\Gamma^\alpha_{\underline{j\alpha}\overline;n}-
    (N+1)(\underset{(1)}{\overset\star{\overline\eta}}{}_{jn}+
    \overline\mu\overline\sigma_{jn})\big)\\&-
    \frac1{N+1}\delta^i_n\big(\overline\Gamma^\alpha_{\underline{j\alpha}\overline;m}-
    (N+1)(\underset{(1)}{\overset\star{\overline\eta}}{}_{jm}+
    \overline\mu\overline\sigma_{jm})\big)\\&
    +
    \big(\overline\sigma_{jm\overline;n}-
    \overline\sigma_{jn\overline;m}-(\overline\sigma_{jm}\overline\sigma_{\alpha n}-
    \overline\sigma_{jn}\overline\sigma_{\alpha m})
    \overline\varphi^\alpha\big)\overline\varphi^i-
    \overline\sigma_{jm}\overline\Gamma^i_{\underset\vee{\alpha
    n}}\overline\varphi^\alpha+
    \overline\sigma_{jn}\overline\Gamma^i_{\underset\vee{\alpha m}}
    \overline\varphi^\alpha
    \\&+\delta^i_j\underset{(1)}{\overset\star\eta}{}_{[mn]}
    -\frac1{N+1}\delta^i_m\big(
    \Gamma^\alpha_{\underline{j\alpha};n}-(N+1)(\underset{(1)}{\overset\star\eta}{}_{jn}+
    \mu\sigma_{jn})\big)\\&+
    \frac1{N+1}\delta^i_n\big(\Gamma^\alpha_{\underline{j\alpha};m}-
    (N+1)(\underset{(1)}{\overset\star\eta}{}_{jm}+\mu\sigma_{jm})\big)\\&
    -
    \big(\sigma_{jm;n}-\sigma_{jn;m}-(\sigma_{jm}\sigma_{\alpha n}-
    \sigma_{jn}\sigma_{\alpha m})\varphi^\alpha\big)\varphi^i+\sigma_{jm}\Gamma^i_{\underset\vee{\alpha
    n}}\varphi^\alpha-
    \sigma_{jn}\Gamma^i_{\underset\vee{\alpha m}}\varphi^\alpha,
    \endaligned\tag{3.1.1}\label{eq:Rtransformationpi31}\\\displaybreak[0]
    &\aligned
    \overline R^i_{jmn}&=R^i_{jmn}-
    \delta^i_j\underset{(2)}{\overset\star{\overline\eta}}{}_{[mn]}
    +\frac1{N+1}\delta^i_m\big(
    \overline\Gamma^\alpha_{\underline{j\alpha}\overline;n}-
    (N+1)(\underset{(2)}{\overset\star{\overline\eta}}{}_{jn}+
    \overline\mu\overline\sigma_{jn})\big)\\&-
    \frac1{N+1}\delta^i_n\big(\overline\Gamma^\alpha_{\underline{j\alpha}\overline;m}-
    (N+1)(\underset{(2)}{\overset\star{\overline\eta}}{}_{jm}+
    \overline\mu\overline\sigma_{jm})\big)\\&
    +
    \big(\overline\sigma_{jm\overline;n}-
    \overline\sigma_{jn\overline;m}-(\overline\sigma_{jm}\overline\sigma_{\alpha n}-
    \overline\sigma_{jn}\overline\sigma_{\alpha m})
    \overline\varphi^\alpha\big)\overline\varphi^i+
    \overline\sigma_{jm}\overline\Gamma^i_{\underset\vee{\alpha
    n}}\overline\varphi^\alpha-
    \overline\sigma_{jn}\overline\Gamma^i_{\underset\vee{\alpha m}}
    \overline\varphi^\alpha
    \\&+\delta^i_j\underset{(2)}{\overset\star\eta}{}_{[mn]}
    -\frac1{N+1}\delta^i_m\big(
    \Gamma^\alpha_{\underline{j\alpha};n}-(N+1)(\underset{(2)}{\overset\star\eta}{}_{jn}+
    \mu\sigma_{jn})\big)\\&+
    \frac1{N+1}\delta^i_n\big(\Gamma^\alpha_{\underline{j\alpha};m}-
    (N+1)(\underset{(2)}{\overset\star\eta}{}_{jm}+\mu\sigma_{jm})\big)\\&
    -
    \big(\sigma_{jm;n}-\sigma_{jn;m}-(\sigma_{jm}\sigma_{\alpha n}-
    \sigma_{jn}\sigma_{\alpha m})\varphi^\alpha\big)\varphi^i-\sigma_{jm}\Gamma^i_{\underset\vee{\alpha
    n}}\varphi^\alpha+
    \sigma_{jn}\Gamma^i_{\underset\vee{\alpha m}}\varphi^\alpha.
    \endaligned\tag{3.1.2}\label{eq:Rtransformationpi32}
  \end{align}

  Hence, based on these transformations and the equation (\ref{eq:T;toT;pi3}) we
  establish the following equations:

  \begin{align}
    &\aligned
    \overline K^i_{jmn}&=K^i_{jmn}-
    \delta^i_j\underset{(1)}{\overset\star{\overline\eta}}{}_{[mn]}
    +\frac1{N+1}\delta^i_m\big(
    \overline\Gamma^\alpha_{\underline{j\alpha}\overline;n}-
    (N+1)(\underset{(1)}{\overset\star{\overline\eta}}{}_{jn}+
    \overline\mu\overline\sigma_{jn})\big)\\&-
    \frac1{N+1}\delta^i_n\big(\overline\Gamma^\alpha_{\underline{j\alpha}\overline;m}-
    (N+1)(\underset{(1)}{\overset\star{\overline\eta}}{}_{jm}+
    \overline\mu\overline\sigma_{jm})\big)\\&
    +
    \big(\overline\sigma_{jm\overline;n}-
    \overline\sigma_{jn\overline;m}-(\overline\sigma_{jm}\overline\sigma_{\alpha n}-
    \overline\sigma_{jn}\overline\sigma_{\alpha m})
    \overline\varphi^\alpha\big)\overline\varphi^i-
    \overline\sigma_{jm}\overline\Gamma^i_{\underset\vee{\alpha
    n}}\overline\varphi^\alpha+
    \overline\sigma_{jn}\overline\Gamma^i_{\underset\vee{\alpha m}}
    \overline\varphi^\alpha
    \\&+\delta^i_j\underset{(1)}{\overset\star\eta}{}_{[mn]}
    -\frac1{N+1}\delta^i_m\big(
    \Gamma^\alpha_{\underline{j\alpha};n}-(N+1)(\underset{(1)}{\overset\star\eta}{}_{jn}+
    \mu\sigma_{jn})\big)\\&+
    \frac1{N+1}\delta^i_n\big(\Gamma^\alpha_{\underline{j\alpha};m}-
    (N+1)(\underset{(1)}{\overset\star\eta}{}_{jm}+\mu\sigma_{jm})\big)\\&
    -
    \big(\sigma_{jm;n}-\sigma_{jn;m}-(\sigma_{jm}\sigma_{\alpha n}-
    \sigma_{jn}\sigma_{\alpha m})\varphi^\alpha\big)\varphi^i+\sigma_{jm}\Gamma^i_{\underset\vee{\alpha
    n}}\varphi^\alpha-
    \sigma_{jn}\Gamma^i_{\underset\vee{\alpha m}}\varphi^\alpha\\&+
    u\underset{(p)}{\overline\sigma}{}^i_{jmn}+
    u'\underset{(q)}{\overline\sigma}{}^i_{jnm}-
    u\underset{(p)}\sigma{}^i_{jmn}-
    u'\underset{(q)}\sigma{}^i_{jnm},
    \endaligned\label{eq:Ktransformationpi31}\\\displaybreak[0]
    &\aligned
    \overline K^i_{jmn}&=K^i_{jmn}-
    \delta^i_j\underset{(2)}{\overset\star{\overline\eta}}{}_{[mn]}
    +\frac1{N+1}\delta^i_m\big(
    \overline\Gamma^\alpha_{\underline{j\alpha}\overline;n}-
    (N+1)(\underset{(2)}{\overset\star{\overline\eta}}{}_{jn}+
    \overline\mu\overline\sigma_{jn})\big)\\&-
    \frac1{N+1}\delta^i_n\big(\overline\Gamma^\alpha_{\underline{j\alpha}\overline;m}-
    (N+1)(\underset{(2)}{\overset\star{\overline\eta}}{}_{jm}+
    \overline\mu\overline\sigma_{jm})\big)\\&
    +
    \big(\overline\sigma_{jm\overline;n}-
    \overline\sigma_{jn\overline;m}-(\overline\sigma_{jm}\overline\sigma_{\alpha n}-
    \overline\sigma_{jn}\overline\sigma_{\alpha m})
    \overline\varphi^\alpha\big)\overline\varphi^i+
    \overline\sigma_{jm}\overline\Gamma^i_{\underset\vee{\alpha
    n}}\overline\varphi^\alpha-
    \overline\sigma_{jn}\overline\Gamma^i_{\underset\vee{\alpha m}}
    \overline\varphi^\alpha
    \\&+\delta^i_j\underset{(2)}{\overset\star\eta}{}_{[mn]}
    -\frac1{N+1}\delta^i_m\big(
    \Gamma^\alpha_{\underline{j\alpha};n}-(N+1)(\underset{(2)}{\overset\star\eta}{}_{jn}+
    \mu\sigma_{jn})\big)\\&+
    \frac1{N+1}\delta^i_n\big(\Gamma^\alpha_{\underline{j\alpha};m}-
    (N+1)(\underset{(2)}{\overset\star\eta}{}_{jm}+\mu\sigma_{jm})\big)\\&
    -
    \big(\sigma_{jm;n}-\sigma_{jn;m}-(\sigma_{jm}\sigma_{\alpha n}-
    \sigma_{jn}\sigma_{\alpha m})\varphi^\alpha\big)\varphi^i-\sigma_{jm}\Gamma^i_{\underset\vee{\alpha
    n}}\varphi^\alpha+
    \sigma_{jn}\Gamma^i_{\underset\vee{\alpha m}}\varphi^\alpha\\&+
    u\underset{(p)}{\overline\sigma}{}^i_{jmn}+
    u'\underset{(q)}{\overline\sigma}{}^i_{jnm}-
    u\underset{(p)}\sigma{}^i_{jmn}-
    u'\underset{(q)}\sigma{}^i_{jnm}.
    \endaligned\label{eq:Ktransformationpi32}
  \end{align}

  \noindent Based on these transformations, we obtain that is

  \begin{eqnarray*}
  \underset{\underset{(1)}{(p,q)}}{\overset\star{\overline{\mathcal
  W}}}{}^i_{jmn}=\underset{\underset{(1)}{(p,q)}}{\overset\star{{\mathcal
  W}}}{}^i_{jmn}&\mbox{and}&
  \underset{\underset{(2)}{(p,q)}}{\overset\star{\overline{\mathcal
  W}}}{}^i_{jmn}=\underset{\underset{(2)}{(p,q)}}{\overset\star{{\mathcal
  W}}}{}^i_{jmn},
  \end{eqnarray*}

  \noindent for $(p,q)\in\{1,\ldots,8\}^2$ and

  \begin{align}
    &\aligned
    \underset{\underset{(1)}{(p,q)}}{\overset\star{{\mathcal
  W}}}{}^i_{jmn}&=K^i_{jmn}+\delta^i_j\underset{(1)}{\overset\star\eta}{}_{[mn]}
    -\frac1{N+1}\delta^i_m\big(
    \Gamma^\alpha_{\underline{j\alpha};n}-(N+1)(\underset{(1)}{\overset\star\eta}{}_{jn}+
    \mu\sigma_{jn})\big)\\&+
    \frac1{N+1}\delta^i_n\big(\Gamma^\alpha_{\underline{j\alpha};m}-
    (N+1)(\underset{(1)}{\overset\star\eta}{}_{jm}+\mu\sigma_{jm})\big)\\&
    -
    \big(\sigma_{jm;n}-\sigma_{jn;m}-(\sigma_{jm}\sigma_{\alpha n}-
    \sigma_{jn}\sigma_{\alpha m})\varphi^\alpha\big)\varphi^i+\sigma_{jm}\Gamma^i_{\underset\vee{\alpha
    n}}\varphi^\alpha-
    \sigma_{jn}\Gamma^i_{\underset\vee{\alpha m}}\varphi^\alpha\\&-
    u\underset{(p)}\sigma{}^i_{jmn}-
    u'\underset{(q)}\sigma{}^i_{jnm},
    \endaligned\label{eq:pi3invpq1}
    \end{align}
    \begin{align}
    &\aligned
    \underset{\underset{(2)}{(p,q)}}{\overset\star{{\mathcal
  W}}}{}^i_{jmn}&=K^i_{jmn}+\delta^i_j\underset{(2)}{\overset\star\eta}{}_{[mn]}
    -\frac1{N+1}\delta^i_m\big(
    \Gamma^\alpha_{\underline{j\alpha};n}-(N+1)(\underset{(2)}{\overset\star\eta}{}_{jn}+
    \mu\sigma_{jn})\big)\\&+
    \frac1{N+1}\delta^i_n\big(\Gamma^\alpha_{\underline{j\alpha};m}-
    (N+1)(\underset{(2)}{\overset\star\eta}{}_{jm}+\mu\sigma_{jm})\big)\\&
    -
    \big(\sigma_{jm;n}-\sigma_{jn;m}-(\sigma_{jm}\sigma_{\alpha n}-
    \sigma_{jn}\sigma_{\alpha m})\varphi^\alpha\big)\varphi^i-\sigma_{jm}\Gamma^i_{\underset\vee{\alpha
    n}}\varphi^\alpha+
    \sigma_{jn}\Gamma^i_{\underset\vee{\alpha m}}\varphi^\alpha\\&-
    u\underset{(p)}\sigma{}^i_{jmn}-
    u'\underset{(q)}\sigma{}^i_{jnm}.
    \endaligned\label{eq:pi3invpq2}
  \end{align}

  \begin{thm}
    Let $f:\mathbb{GR}_N\rightarrow\mathbb{G\overline R}_N$ be an
    equitorsion almost geodesic mapping of the third type and $s$-th
    kind, $s=1,2$, which has the property of reciprocity. Families
    $\underset{\underset{(1)}{(p,q)}}{\overset\star{{\mathcal
        W}}}{}^i_{jmn}$ and
        $\underset{\underset{(2)}{(p,q)}}{\overset\star{{\mathcal
        W}}}{}^i_{jmn}$ given by the equations \emph{(\ref{eq:pi3invpq1},
        \ref{eq:pi3invpq2})} are families of invariants of mapping of the
        corresponding kind.\qed
  \end{thm}

  \begin{cor}
    The families $\underset{\underset{(s)}{(p,q)}}{\overset\star{{\mathcal
        W}}}{}^i_{jmn}, s=1,2$,  of invariants of an equitorsion almost
        geodesic mapping
        $f:\mathbb{GR}_N\rightarrow\mathbb{G\overline R}_N$ which
        has the property of reciprocity and the invariants
        $\underset{(s)}{\overset\star{\mathcal W}}{}^i_{jmn}$ given
        by the equations \emph{(\ref{eq:inv1}, \ref{eq:inv2})}
        satisfy the equations

        \begin{equation}
        \aligned
          \underset{\underset{s}{(p,q)}}{\overset\star{\mathcal
          W}}{}^i_{jmn}&=\underset{(s)}{\overset\star{\mathcal
          W}}{}^i_{jmn}-u\underset{(p)}\sigma{}^i_{jmn}-
          u'\underset{(q)}\sigma{}^i_{jnm}\\&+u\Gamma^i_{\underset\vee{jm};n}+
          u'\Gamma^i_{\underset\vee{jn};m}+
          v\Gamma^\alpha_{\underset\vee{jm}}\Gamma^i_{\underset\vee{\alpha
          n}}+
          v'\Gamma^\alpha_{\underset\vee{jn}}\Gamma^i_{\underset\vee{\alpha
          m}}+
          w\Gamma^\alpha_{\underset\vee{mn}}\Gamma^i_{\underset\vee{\alpha
          j}},
        \endaligned\label{eq:pi3Wcorrelations}
        \end{equation}

        \noindent for $(p,q)\in\{1,\ldots,8\}^2$.\qed
  \end{cor}
  \begin{cor}
    The rank of matrix

    \begin{equation}
      \overset\star{\mathcal W}=\left[\begin{array}{ccccccccc}
        1&-u^\rho_1&\ldots&-u^\rho_{20}&u&u'&v&v'&w
      \end{array}
      \right]
    \end{equation}

    \noindent of the type $64\times26$ is equal $6$, i.e. there are
    six linearly independent families
    $\underset{\underset{(s)}{(p,q)}}{\overset\star{\mathcal
    W}}{}^i_{jmn},s\in\{1,2\},\linebreak(p,q)\in\{1,\ldots,8\}^2$, of invariants
    given by the equations \emph{(\ref{eq:pi3invpq1},
    \ref{eq:pi3invpq2})}.\qed
  \end{cor}

\end{document}